\RequirePackage{etoolbox}

\documentclass[ba,preprint]{imsart}
\pubyear{2017}
\volume{12}
\issue{4}
\doi{10.1214/17-BA1065}
\firstpage{1262}
\lastpage{1263}

\usepackage{amsthm}
\usepackage{amsmath}
\usepackage{natbib}
\usepackage[colorlinks,citecolor=blue,urlcolor=blue,filecolor=blue,backref=page]{hyperref}
\usepackage{graphicx}

\startlocaldefs
\endlocaldefs

\begin{document}


\begin{frontmatter}
\title{Contributed Discussion to Uncertainty Quantification for the Horseshoe by St\'{e}phanie van der Pas, Botond Szab\'{o} and Aad van der Vaart}

\runtitle{Horseshoe Discussion}

\begin{aug}
\author{\fnms{William Weimin} \snm{Yoo}\thanksref{addr1}\ead[label=e1]{yooweimin0203@gmail.com}}

\runauthor{W. W. Yoo}

\address[addr1]{Mathematical Institute, Leiden University, The Netherlands \printead{e1}}


\end{aug}

\begin{abstract}
We begin by introducing the main ideas of the paper under discussion. We discuss some interesting issues regarding adaptive component-wise credible intervals. We then briefly touch upon the concepts of self-similarity and excessive bias restriction. This is then followed by some comments on the extensive simulation study carried out in the paper.
\end{abstract}

\begin{keyword}
\kwd{Horseshoe}
\kwd{Half-Cauchy}
\kwd{MMLE}
\kwd{Credible intervals}
\kwd{Model Selection}
\kwd{Credible balls}
\kwd{Adaptive}
\kwd{Excessive bias restriction}
\end{keyword}

\end{frontmatter}

I would like to congratulate the authors for such an comprehensive and interesting paper on the horseshoe prior and its use in Bayesian uncertainty quantification. Let me first summarize key ideas of \cite{horseshoe} in order to set the tone for my discussion. The horseshoe prior is a scale mixture of normals with the half-Cauchy as the mixing distribution. The half-Cauchy is in turn the absolute value of a standard Cauchy distribution. Working under the normal means model $Y_i=\theta_i+\varepsilon_i,i=1,\dotsc,n$, this hierarchical prior takes the form $\theta_i|\lambda_i,\tau\sim\mathrm{N}(0,\lambda_i^2\tau^2)$ and $\lambda_i\sim\text{Half-Cauchy}$. The paper considers two methods of estimating $\tau$, namely by empirical Bayes through the maximum marginal likelihood estimator (MMLE) or by endowing another layer of hyper-prior.

The true signal $\theta_0$ is assumed to be sparse, and recovery of the nonzero values using the horseshoe prior was studied by the same authors previously in \cite{horseshoerate}. The present paper however deals with the issue of accessing the quality of this recovery procedure. This is accomplished through the construction of (adaptive) component-wise credible intervals and $\ell_2$-credible balls. Moreover, the authors introduced a simple model selection procedure by declaring that a signal component is unimportant if its corresponding credible interval includes $0$ within its span.

My discussion will focus on the case of component-wise credible interval, since I find that its results are the most interesting and these intervals are the ones used in the simulations. The most prominent feature regarding results for these intervals is the division of the true signal components into three regimes corresponding to small, intermediate and large signals. The horseshoe interval was able to provide adequate coverage for small and large signals but not the intermediate ones. The existence of this intermediate layer and the gaps in between these regimes made me wonder whether this is due to the intrinsic nature of component-wise credible intervals, or other more extrinsic factors such as the horseshoe prior used or perhaps an artefact of the proof techniques employed.

It is now well known that adaptive credible sets cannot do honest uncertainty quantification over all possible true signals, and some of these signals must be permanently excluded. To this end, the authors discussed two criteria for removal, one based on the concept of self-similarity and the other based on the excessive bias restriction introduced by \cite{ebr}. In the present setting of sparse signals, the key insight into these conditions is that the true signals must be at some distance away from the zero signal in a sense made precise in the paper. Interestingly as mentioned in Remark 3, the three regimes become more ``contiguous" under self-similarity when compared to the situation where self-similarity was not assumed, as this is evident by comparing $\mathcal{S}_a,\mathcal{M}_a,\mathcal{L}_a$ with $\mathcal{S},\mathcal{M},\mathcal{L}$ for $\tau=\tau_n(p_n)$. This in turn suggests that throwing away troublesome truths enables the horseshoe credible interval to fill in the gaps between the three regions.

For the sake of discussion, let us continue working under the self-similarity or excessive bias restriction. From the simulation results, it is clear that the horseshoe credible intervals have the best performance in terms of high coverage and shortest lengths when the means (or the true signals) are zero. Two settings of $p_n$ the number of nonzero signals were used, i.e., $p_n=20$ and $p_n=200$ when $n=400$. Now let us increase the proportion of zero means ($n-p_n$) to a point that the self-similarity condition is violated, will the horseshoe still enjoy this near-perfect performance? By looking at the bar charts on coverage and lengths, it is conceivable that we can still get good performance even if this condition is violated slightly. My question concerns whether it is possible to observe empirically what will happen when sparse signals become non self-similar in the sense discussed in the paper.

Uncertainty quantification is undoubtedly one of the most active research areas in Bayesian statistics, and as the present paper shows, it involves resolving many delicate technical and practical issues. The horseshoe prior has proven itself to be optimal in sparse signal recovery, and we are able to access the quality of this recovery thanks to the theories and methods developed in the paper. It would be interesting also to consider other classes of priors, e.g., spike-and-slab types, and I hope that there will be more papers on Bayesian uncertainty quantification for sparse models in the future.




\bibliographystyle{ba}
\bibliography{horseshoeref}


\end{document}